\documentclass[reqno,12pt]{amsart}
\usepackage{amssymb}
\usepackage{amsmath}
\usepackage{amsfonts}

\newtheorem{theorem}{Theorem}[section]
\theoremstyle{plain}

\newtheorem{corollary}[theorem]{Corollary}

\newtheorem{definition}[theorem]{Definition}

\newtheorem{lemma}[theorem]{Lemma}

\newtheorem{remark}[theorem]{Remak}

\numberwithin{equation}{section}
\input{tcilatex}

\begin{document}
\title{Convexity in the interpolation spaces }
\author{Daher Mohammad}
\address{117, Rue Gustave Courbet 77350 Le M\'{e}e Sur Seine-France}
\email{daher.mohammad@ymail.com}

\begin{abstract}
In this work we study if the norms rotund, uniformly rotund, weakly
uniformly rotund, locally uniformly rotund or weakly locally uniformly
rotund interpolate in the complex or the real interpolation spaces. We will
see that the properties uniformly rotund and weakly uniformly rotund
interpolate, but the properties locally uniformly rotund or weakly locally
unifomly rotund do not interpolate in general except for the dual
interpolation couple.
\end{abstract}

\maketitle

\section{\protect\bigskip \textsc{Introduction.}\ \ \ \ \ \ \ \ \ \ \ \ \ \
\ \ \ \ \ \ \ \ \ \ \ \ \ \ \ \ \ \ \ \ \ \ \ \ \ \ \ \ \ \ \ \ \ \ \ \ \ \
\ \ \ \ \ \ \ \ \ \ \ \ \ \ \ \ \ \ \ \ \ \ \ \ \ \ \ \ \ \ \ \ \ \ \ \ \ \
\ \ \ \ \ \ \ \ \ \ \ \ \ \ \ \ \ \ \ \ \ \ \ \ \ \ \ \ \ \ \ \ \ \ \ \ \ \
\ \ \ \ \ \ \ \ \ \ \ \ \ \ \ \ \ \ \ \ \ \ \ \ \ \ \ \ \ \ \ \ \ \ \ \ \ \
\ \ \ \ \ \ \ \ \ \ \ \ \ \ \ \ \ \ \ \ \ \ \ \ \ \ \ \ \ \ \ \ \ \ \ \ \ \
\ \ \ \ \ \ \ \ \ \ \ \ \ \ \ \ \ \ \ \ \ \ \ \ \ \ \ \ \ \ \ \ \ \ \ \ \ \
\ \ \ \ \ \ \ \ \ \ \ \ \ \ \ \ \ \ \ \ \ \ \ \ \ \ \ \ \ \ \ \ \ \ \ \ \ \
\ \ \ \ \ \ \ \ \ \ \ \ \ \ \ \ \ \ \ \ \ \ \ \ \ \ \ \ \ \ \ \ \ \ \ \ \ \
\ \ \ \ \ \ \ \ \ \ \ \ \ \ \ \ \ \ \ \ \ \ \ \ \ \ \ \ \ \ \ \ \ \ \ \ \ \
\ \ \ \ \ \ \ \ \ \ \ \ \ \ \ \ \ \ \ \ \ \ \ \ \ \ \ \ \ \ \ \ \ \ \ \ \ \
\ \ \ \ \ \ \ \ \ \ \ \ \ \ \ \ \ \ \ \ \ \ \ \ \ \ \ \ \ \ \ \ \ \ \ \ \ \
\ \ \ \ \ \ \ \ \ \ \ \ \ \ \ \ \ \ \ \ \ \ \ \ \ \ \ \ \ \ \ \ \ \ \ \ \ \
\ \ \ \ \ \ \ \ \ \ \ \ \ \ \ \ \ \ \ \ \ \ \ \ \ \ \ \ \ \ \ \ \ \ \ \ \ \
\ \ \ \ \ \ \ \ \ \ \ \ \ \ \ \ \ \ \ \ \ \ \ \ \ \ \ \ \ \ \ \ \ \ \ \ \ \
\ \ \ \ \ \ \ \ \ \ \ \ \ \ \ \ \ \ \ \ \ \ \ \ \ \ \ \ \ \ \ \ \ \ \ \ \ \
\ \ \ \ \ \ \ \ \ \ \ \ \ \ \ \ \ \ \ \ \ \ \ \ \ \ \ \ \ \ \ \ \ \ \ \ \ \
\ \ \ \ \ \ \ \ \ \ \ \ \ \ \ \ \ \ \ \ \ \ \ \ \ \ \ \ \ \ \ \ \ \ \ \ \ \
\ \ \ \ \ \ \ \ \ \ \ \ \ \ \ \ \ \ \ \ \ \ \ \ \ \ \ \ \ \ \ \ \ \ \ \ \ \
\ \ \ \ \ \ \ \ \ \ \ \ \ \ \ \ \ \ \ \ \ \ \ \ \ \ \ \ \ \ \ \ \ \ \ \ \ \
\ \ \ \ \ \ \ \ \ \ \ \ \ \ \ \ \ \ \ \ \ \ \ \ \ \ \ \ \ \ \ \ \ \ \ \ \ \
\ \ \ \ \ \ \ \ \ \ \ \ \ \ \ \ \ \ \ \ \ \ \ \ \ \ \ \ \ \ \ \ \ \ \ \ \ \
\ \ \ \ \ \ \ \ \ \ \ \ \ \ \ \ \ \ \ \ \ \ \ \ \ \ \ \ \ \ \ \ \ \ \ \ \ \
\ \ \ \ \ \ \ \ \ \ \ \ \ \ \ \ \ \ \ \ \ \ \ \ \ \ \ \ \ \ \ \ \ \ \ \ \ \
\ \ \ \ \ \ \ \ \ \ \ \ \ \ \ \ \ \ \ \ \ \ \ \ \ \ \ \ \ \ \ \ \ \ \ \ \ \
\ \ \ \ \ \ \ \ \ \ \ \ \ \ \ \ \ \ \ \ \ \ \ \ \ \ \ \ \ \ \ \ \ \ \ \ \ \
\ \ \ \ \ \ \ \ \ \ \ \ \ \ \ \ \ \ \ \ \ \ \ \ \ \ \ \ \ \ \ \ \ \ \ \ \ \
\ \ \ \ \ \ \ \ \ \ \ \ \ \ \ \ \ \ \ \ \ \ \ \ \ \ \ \ \ \ \ \ \ \ \ \ \ \
\ \ \ \ \ \ \ \ \ \ \ \ \ \ \ \ \ \ \ \ \ \ \ \ \ \ \ \ \ \ \ \ \ \ \ \ \ \
\ \ \ \ \ \ \ \ \ \ \ \ \ \ \ \ \ \ \ \ \ \ \ \ \ \ \ \ \ \ \ \ \ \ \ \ \ \
\ \ \ \ \ \ \ \ \ \ \ \ \ \ \ \ \ \ \ \ \ \ \ \ \ \ \ \ \ \ \ \ \ \ \ \ \ \
\ \ \ \ \ \ \ \ \ \ \ \ \ \ \ \ \ \ \ \ \ \ \ \ \ \ \ \ \ \ \ \ \ \ \ \ \ \
\ \ \ \ \ \ \ \ \ \ \ \ \ \ \ \ \ \ \ \ \ \ \ \ \ \ \ \ \ \ \ \ \ \ \ \ \ \
\ \ \ \ \ \ \ \ \ \ \ \ \ \ \ \ \ \ \ \ \ \ \ \ \ \ \ \ \ \ \ \ \ \ \ \ \ \
\ \ \ \ \ \ \ \ \ \ \ \ \ \ \ \ \ \ \ \ \ \ \ \ \ \ \ \ \ \ \ \ \ \ \ \ \ \
\ \ \ \ \ \ \ \ \ \ \ \ \ \ \ \ \ \ \ \ \ \ \ \ \ \ \ \ \ \ \ \ \ \ \ \ \ \
\ \ \ \ \ \ \ \ \ \ \ \ \ \ \ \ \ \ \ \ \ \ \ \ \ \ \ \ \ \ \ \ \ \ \ \ \ \
\ \ \ \ \ \ \ \ \ \ \ \ \ \ \ \ \ \ \ \ \ \ \ \ \ \ \ \ \ \ \ \ \ \ \ \ \ \
\ \ \ \ \ \ \ \ \ \ \ \ \ \ \ \ \ \ \ \ \ \ \ \ \ \ \ \ \ \ \ \ \ \ \ \ \ \
\ \ \ \ \ \ \ \ \ \ \ \ \ \ \ \ \ \ \ \ \ \ \ \ \ \ \ \ \ \ \ \ \ \ \ \ \ \
\ \ \ \ \ \ \ \ \ \ \ \ \ \ \ \ \ \ \ \ }

Let $\overline{B}=(B_{0},B_{1})$ be an interpolation couple, let $\theta \in
(0,1)$ and let $p\in (1,+\infty )$. If the complex interpolation space $%
B_{\theta }$ (or the real interpolation space B$_{\theta ,p})$ preserves a
given geometric property of $B_{0}$ or $B_{1},$ then we say that this
property is interpolate. For example, the Reflexivity and the separability
are interpolated. In the parts 4 (resp. part 6), we show that if the norm on 
$B_{0}^{\ast }$ is rotund (resp. $LUR$ or weakly $LUR)$, then the norm on $%
B_{\theta }^{\ast }$ is rotund (resp. is $LUR$ or weakly $LUR),$ 0%
\TEXTsymbol{<}$\theta <1.$ In the part 5 we show that if the norm on $B_{0}$
is $UR,$ then the norm on $B_{\theta }$ is $UR,$ $0<\theta <1.$ In the part
7 and 8 we show similar results concering the real interpolation spaces
(modulo an equivalent norm).

In \cite{Da2}, we show that if the norm on $B_{0}$ is $WUR,$ then the norm
on $B_{\theta }$ is $WUR.$ We recall by \cite{Da3}, there exists an
interpolation couple $(A_{0},A_{1})$ such that the norm on $A_{0}$ is $LUR$
but the interpolated spaces $A_{\theta },A_{\theta ,p}$ not admit any
equivalent rotund norm for all $0<\theta <1$ and all $1<p<+\infty $.

\section{Definitions and some facts about complex interpolation}

We denote by $\left\langle x,x^{\ast }\right\rangle $ the pairing between an
element $x$ of a Banach space $X$ and an element $x^{\ast }$ of its dual $%
X^{\ast }$. Let $S_{0}=\left\{ z\in \mathbb{C}\mid \ 0<\mathop{\rm Re}%
\nolimits(z)<1\right\} $ and let $S$ be its closure. Let $\overline{B}%
=(B_{0},B_{1})$ be a complex interpolation couple, in the sense of \cite[%
Chap.~II, p.~24]{Ber-Lof}. We first recall the definition of the
interpolation space $B_{\theta }$, $\theta \in (0,1)$ \cite[Chap.~4]{Ber-Lof}%
.

Let $\mathcal{F}(\overline{B})$ be the space of functions $F$ with values in 
$B_{0}+B_{1}$, which are bounded and continuous on $S$, holomorphic on $%
S_{0} $, such that, for $j\in \{0,1\}$, the map $F(j + i \, \cdot)$ lies in $%
\mathcal{C}_{0}(\mathbb{R},B_{j})$. We equip $\mathcal{F}(\overline{B})$
with the norm

\begin{equation*}
\bigl\Vert F\bigr\Vert_{\mathcal{F}(\overline{B})}^{{}}=\max_{j\in \left\{
0,1\right\} }\sup_{\tau \in \mathbb{R}}\bigl\Vert F(j+i\tau )\bigr\Vert%
_{B_{j}}^{{}}.
\end{equation*}

\noindent The space $B_{\theta }=(B_{0},B_{1})_{\theta }=\left\{ F(\theta
)\mid \ F\in \mathcal{F}(\overline{B})\right\}$, $0<\theta <1$, is a Banach
space \cite[Th. 4.1.2]{Ber-Lof} for the norm defined by

\begin{equation*}
\bigl\Vert a\bigr\Vert_{B_{\theta }}^{{}}=\inf \left\{ \bigl\Vert F\bigr\Vert%
_{\mathcal{F}(\overline{B})}^{{}};\mathrm{\ }F(\theta )=a\right\} .
\end{equation*}

\bigskip Every $F\in \mathcal{F}\mathfrak{(}\overline{B})$ have the
representation froms its values on the boundary of $S_{0}$ ($\left\{ j+i\tau
;\text{ }j\in \left\{ 0,1\right\} ,\text{ }\tau \in \mathbb{R}\right\} )$
thanks to the harmonic measure \cite[Sections 4.3, 4.5]{Ber-Lof}$:$ if $%
z=\theta +it,$ $\frac{Q_{0}(z,.)}{1-\theta }$ and $\frac{Q_{1}(z,.)}{\theta }
$ are the probability densities on $\mathbb{R}$.

\begin{equation}
F(w)=\mathop{\displaystyle \int}\limits_{\mathbb{R}}F(i\tau )Q_{0}(w,\tau
)\,d\tau +\mathop{\displaystyle \int}\limits_{\mathbb{R}}F(1+i\tau
)Q_{1}(w,\tau )\,d\tau ,w\in S_{0}.  \label{1}
\end{equation}

Let $\mathcal{H}(B_{0},B_{1})$ be the space of functions $F:S\rightarrow
B_{0}+B_{1}$ which are holomorphic on $S_{0}$, such that for $j\in \{0,1\}$, 
$\tau \rightarrow F(j+i\tau )$ is a.s. the weak non-tangential limit of $%
F(z) $ as $z\rightarrow j+i\tau $ and $\tau \rightarrow F(j+i\tau )$ is
strongly measurable with values in $B_{j},$ $j\in \left\{ 0,1\right\} $ (in
the sense of \cite[Chap.~II, p.~41]{DU}).

\bigskip For $1\leq p\leq \infty $ and $\theta \in (0,1)$, $\mathcal{F}%
_{\theta }^{p}(B_{0},B_{1})$ denotes the space of functions $F\in \mathcal{H}%
(B_{0},B_{1})$ such that if $p=+\infty ,$

\begin{equation*}
\bigl\Vert F\bigr\Vert_{\mathcal{F}_{\theta }^{\infty
}(B_{0},B_{1})}^{{}}=\max \{\bigl\Vert F(i.)\bigr\Vert_{L^{\infty
}(Q_{0}(\theta ,.)d\tau ,B_{0})}^{{}},\bigl\Vert F(1+i.)\bigr\Vert%
_{L^{\infty }(Q_{1}(\theta ,.)d\tau ,B_{1})}^{{}}\}<+\infty ,
\end{equation*}

\noindent\ or, if $p<+\infty ,$ 
\begin{equation*}
\bigl\Vert F\bigr\Vert_{\mathcal{F}_{\theta }^{p}(B_{0},B_{1})}^{p}=%
\mathop{\displaystyle \int}\limits_{\mathbb{R}}\bigl\Vert F(i\tau )\bigr\Vert%
_{B_{0}}^{p}Q_{0}(\theta ,\tau )d\tau +\mathop{\displaystyle \int}\limits_{%
\mathbb{R}}\bigl\Vert F(1+i\tau )\bigr\Vert_{B_{1}}^{p}Q_{1}(\theta ,\tau
)d\tau <+\infty ,
\end{equation*}

and $F$ satisfying (\ref{1}). $\mathcal{F}_{\theta }^{p}(B_{0},B_{1})$ deos
not depend of $\theta ,$ because $Q_{0}(\theta ,.)$ and $Q_{1}(\theta ,.)$
are continuous and stricly positives on $\mathbb{R}.$ As $\mathcal{F}%
_{\theta }^{p}(B_{0},B_{1})=\mathcal{F}_{\theta ^{\prime }}^{p}(B_{0},B_{1})$
isometrically, we will denote sometimes $\mathcal{F}_{\theta }^{\infty
}(B_{0},B_{1})=\mathcal{F}^{\infty }(B_{0},B_{1}).$ $\mathcal{F}%
(B_{0},B_{1}) $ is an isometrically subspace of $\mathcal{F}^{\infty
}(B_{0},B_{1}).$ \bigskip

Let $(B_{0},B_{1})$ be an interpolation couple. The dual of $B_{j},$ $j\in
\left\{ 0,1\right\} $ (resp. $B_{\theta })$ will be denoted by $B_{j}^{\ast
},$ $j\in \left\{ 0,1\right\} $ (resp. $B_{\theta }^{\ast })$, $0<\theta <1.$

We now recall the definition of the ``upper" interpolation space $B^{\theta
} $. Let $\mathcal{G(}\overline{B})$ be the space of functions $%
G:S\rightarrow B_{0}+B_{1}$, which are continuous on $S$, holomorphic on $%
S_{0}$, such that

(i) the map: $z\rightarrow (1+\left\vert z\right\vert )^{-1}\left\Vert
G(z)\right\Vert _{B_{0}+B_{1}}$ is bounded on $S$.

(ii) $G(j+i\tau )-G(j+i\tau ^{\prime })\in B_{j}$ for every $\tau ,\tau
^{\prime }\in \mathbb{R}$, $j\in \{0,1\}$, and 
\begin{equation*}
\max_{j\in \left\{ 0,1\right\} }\sup_{\tau \neq \tau ^{\prime }\in \mathbb{R}%
}\biggl(\frac{\bigl\Vert G(j+i\tau )-G(j+i\tau ^{\prime })\bigr\Vert%
_{B_{j}}^{{}}}{\left\vert \tau -\tau ^{\prime }\right\vert }\biggr)<+\infty .
\end{equation*}

\noindent By \cite[Lemma 4.1.3]{Ber-Lof} this defines a norm $\bigl\Vert\dot{%
G}\bigr\Vert_{Q\mathcal{G}(\overline{B})}^{{}}$ on $Q\mathcal{G(}\overline{B}%
)$, the quotient of $\mathcal{G}(\overline{B})$ by the subspace of constant
functions (in $B_{0}+B_{1}$), and $Q\mathcal{G}(\overline{B})$ is complete
with respect to this norm.

The space $B^{\theta }=\left\{ a\in B_{0}+B_{1}\mid \ \exists \;G\in 
\mathcal{G(}\overline{B}),\ a=G^{\prime }(\theta )\right\} $ is a Banach
space \cite[Th. 4.1.4]{Ber-Lof} with respect to the norm

\begin{equation*}
\bigl\Vert a\bigr\Vert_{B^{\theta }}^{p}=\inf \left\{ \bigl\Vert\dot{G}%
\bigr\Vert_{Q\mathcal{G}(\overline{B})}^{{}}\mid \ G^{\prime }(\theta
)=a\right\} .
\end{equation*}

\noindent We shall use the following properties (and will not explicitly
quote them):

1) $B_{0}\cap B_{1}$ is dense in $B_{\theta }$ \cite[Th.~4.2.2]{Ber-Lof}.

\noindent We recall that $(B_{0},B_{1})$ is a regular couple if $B_{0}\cap
B_{1}$ is dense in $B_{0}$ and $B_{1}$. Then, for a regular couple,

2) the dual of $B_{0}\cap B_{1}$ is $B_{0}^{\ast }+B_{1}^{\ast }$
isometrically \cite[Th.~2.7.1]{Ber-Lof}.

3) the dual of $B_{\theta }$ is $(B_{0}^{\ast },B_{1}^{\ast })^{\theta }$
isometrically \cite[Th.~4.5.1]{Ber-Lof}. \bigskip

\section{\textsc{preliminary}}

\bigskip We now assume that $(B_{0},B_{1})$ is a regular couple (hence ($%
B_{0}\cap B_{1})^{\ast }=B_{0}^{\ast }+B_{1}^{\ast }$) and that $B_{0},B_{1}$
are separable spaces.

\noindent We recall that $B_{0}\cap B_{1}$ is then separable. Indeed the map 
$\sigma :(x,y)\rightarrow x+y$, $B_{0}^{\ast }\times B_{1}^{\ast
}\rightarrow B_{0}^{\ast }+B_{1}^{\ast }$ is continuous for the respective $%
w^{\ast }$-topologies. The norm closed unit balls $X,Y$ of these spaces,
equipped with the $w^{\ast }$-topologies, are compact and $X$ is metrizable.
The set $K=X\cap \sigma ^{-1}(Y)$ is closed in $X$, hence compact metrizable
for the $w^{\ast }$-topology of $X$. Since $\sigma (K)=Y$ \ by definition of
the norm in $B_{0}^{\ast }+B_{1}^{\ast },$ the map $f\rightarrow f\circ
\sigma $ identifies isometrically $\mathcal{C}(Y)$ with a norm closed
subspace of the separable space $\mathcal{C}(K)$. Hence $B_{0}\cap B_{1},$
being a norm closed subspace of $\mathcal{C}(Y),$ is separable.

\noindent The definition of $\mathcal{F}_{\ast }^{\infty }(B_{0}^{\ast
},B_{1}^{\ast })$ will use the following known fact, for which we prefer to
precise a proof.

\emph{Fact : }Let $B$ be a separable Banach space and let $%
g:S_{0}\rightarrow B^{\ast }$ be a bounded holomorphic function. Then $g$
admits an a.s. defined extension $g_{\ast }$ on the boundary of $S_{0},$
with values in $B^{\ast }$, and $g_{\ast }$ is $\sigma (B^{\ast },B)$
measurable.

Indeed, for every $a\in B$, the function $w\rightarrow \left\langle
a,g(w)\right\rangle $ has, for almost every $\tau \in \mathbb{R}$, a non
tangential limit as $w\rightarrow i\tau ,$ denoted by $l_{a}(i\tau )$, and $%
\bigl\Vert l_{a}\bigr\Vert_{L^{\infty }(id\tau )}^{{}}\leq \bigl\Vert g%
\bigr\Vert_{\infty }^{{}}\bigl\Vert a\bigr\Vert_{B}^{{}}$. For $a,b\in B,$
if $l_{a}(i\tau )$, $l_{b}(i\tau )$ are defined, then $l_{\lambda a+\mu
b}(i\tau )=\lambda l_{a}(i\tau )+\mu l_{b}(i\tau )$ for every scalars $%
\lambda ,\mu .$ Hence, if $b_{n}\rightarrow _{n\rightarrow \infty }a$ in $B,$
$\bigl\Vert l_{b_{n}}-l_{a}\bigr\Vert_{L^{\infty }(id\tau )}^{{}}\rightarrow
_{n\rightarrow \infty }0$. Since\ $B$ is separable let $(a_{n})_{n\geq 1}$
be a dense sequence in $B,$ let $V$ be their linear span$.$ For almost every 
$\tau ,$ the map : $a_{n\rightarrow }l_{a_{n}}(i\tau )$ extends as a linear
form on $V$ with norm less than $\left\Vert g\right\Vert _{\infty },$ hence
defines an element $g_{\ast }(i\tau )$ in $B^{\ast }$ with $\bigl\Vert %
g_{\ast }(i\tau )\bigr\Vert_{B^{\ast }}^{{}}\leq \bigl\Vert g\bigr\Vert%
_{\infty }^{{}}$, such that $l_{b}(i\tau )=\left\langle b,g_{\ast }(i\tau
)\right\rangle $ for every $b\in V;$ finally, by density of the $a_{n}$'s in 
$B,$ such an equality holds for every $b\in B.$ The proof is similar for $1+i%
\mathbb{R}$.

Taking $B=B_{0}\cap B_{1}$ we denote by $\mathcal{F}_{\ast }^{\infty
}(B_{0}^{\ast },B_{1}^{\ast })$ the space of these functions $g$ such that
moreover a.s. $\tau \rightarrow g_{\ast }(j+i\tau )$ is bounded : $\mathbb{R}%
\rightarrow B_{j}^{\ast }$, $j\in \left\{ 0,1\right\} .$

\noindent In particular, $g_{\ast }(j+i\,\cdot )$ is $\sigma (B_{j}^{\ast
},B_{j})$ measurable since $B_{0}\cap B_{1}$ is dense in $B_{j}$ and the
function $\tau \rightarrow \bigl\Vert g_{\ast }(j+i\,\cdot )\bigr\Vert%
_{B_{j}^{\ast }}^{{}}$ is measurable since $B_{j}$ is separable, $\ j\in
\left\{ 0,1\right\} $. This allows to equip the space $\mathcal{F}_{\ast
}^{\infty }(B_{0}^{\ast },B_{1}^{\ast })$ with the norm 
\begin{equation*}
\bigl\Vert g\bigr\Vert_{\mathcal{F}_{\ast }^{\infty }(B_{0}^{\ast
},B_{1}^{\ast })}^{{}}=\max_{j\in \left\{ 0,1\right\} }\bigl\Vert\bigl\Vert %
g_{\ast }(j+i\,\cdot )\bigr\Vert_{_{B_{j}^{\ast }}}\bigr\Vert_{L^{\infty
}(d\tau )}^{{}}.
\end{equation*}

\noindent Since for every $a\in B_{0}\cap B_{1}$, $\left\langle
a,g(.)\right\rangle $ satisfies the integral representation (\ref{1}), we
get 
\begin{equation*}
\bigl\Vert g(w)\bigr\Vert_{B_{0}^{\ast }+B_{1}^{\ast }}\leq \bigl\Vert g%
\bigr\Vert_{\mathcal{F}_{\ast }^{\infty }(B_{0}^{\ast },B_{1}^{\ast
})},\quad w\in S_{0}.
\end{equation*}

\begin{lemma}
\label{JM}\cite[Lemma 3.2]{Da5}Let $(B_{0},B_{1})$ be a regular couple of
separable spaces, let $\theta \in S_{0}$ and let $a^{\ast }\in (B_{0}^{\ast
},B_{1}^{\ast })^{\theta }$. Then%
\begin{eqnarray*}
\left\Vert a^{\ast }\right\Vert _{B_{\theta }^{\ast }} &=&\inf \left\{
\left\Vert h\right\Vert _{\mathcal{F}_{\ast }^{\infty }(B_{0}^{\ast
},B_{1}^{\ast })};\ h(\theta )=a^{\ast }\right\} \\
&=&\inf \biggl\{\underset{j=0}{\overset{1}{\mathop{\displaystyle \sum }}}%
\mathop{\displaystyle \int}\limits_{\mathbb{R}}\bigl\Vert h_{\ast }(i\tau )%
\bigr\Vert_{_{B_{j}^{\ast }}}Q_{j}(\theta ,\tau )d\tau ;\ h\in \mathcal{F}%
_{\ast }^{\infty }(B_{0}^{\ast },B_{1}^{\ast }),\ h(\theta )=a^{\ast }%
\biggr\}.
\end{eqnarray*}
\end{lemma}

\begin{lemma}
\label{YO}\cite[Lemma 3.3]{Da5}Let $(B_{0},B_{1})$ be a regular couple of
separable spaces, let $\theta \in (0;1)$ and let $a^{\ast }\in (B_{0}^{\ast
},B_{1}^{\ast })^{\theta }$ such that $a^{\ast }\neq 0$. Then there exists $%
h\in \mathcal{F}_{\ast }^{\infty }(B_{0}^{\ast },B_{1}^{\ast }),\ ^{\ast }$
such that $h(\theta )=a$ and $\left\Vert a^{\ast }\right\Vert _{(B_{0}^{\ast
},B_{1}^{\ast })^{\theta }}=\underset{j=0}{\overset{1}{%
\mathop{\displaystyle
\sum }}}\mathop{\displaystyle \int}\limits_{\mathbb{R}}\bigl\Vert h_{\ast
}(i\tau )\bigr\Vert_{_{B_{j}^{\ast }}}Q_{j}(\theta ,\tau )d\tau $.
\end{lemma}

Let $(B_{0},B_{1})$ be an interpolation couple and let $r\geq 1.$ On $%
B_{0}+B_{1}$ one defines the norm $\left\vert .\right\vert _{r}$ by $%
\left\vert a\right\vert _{r}=\inf \left\{ (\left\Vert a_{0}\right\Vert
_{B_{0}}^{r}+\left\Vert a_{1}\right\Vert _{B_{1}}^{r})^{\frac{1}{r}};\text{ }%
a_{0}+a_{1}=a\right\} ,$ $a\in B_{0}+B_{1}.$ One defines also the norm $%
N_{r}(.)$ on $B_{0}\cap B_{1}$ by $N_{r}(a)=(\left\Vert a\right\Vert
_{B_{0}}^{r}+\left\Vert a\right\Vert _{B_{1}}^{r})^{\frac{1}{r}},$ $a\in
B_{0}\cap B_{1}.$

For $r\geq 1,$ we denote by $r^{\prime }$ the conjugate of $r.$

\begin{lemma}
\label{Dual}\cite[Lemme 2.1]{Da6}Let $(B_{0}^{{}},B_{1}^{{}})$ be a regular
couple and let $r\geq 1$. Then

1) $(B_{0}+B_{1},\left\vert .\right\vert _{r})^{\ast }=(B_{0}^{\ast }\cap
B_{1}^{\ast },N_{r^{\prime }}(.))$ isometrically.

2) $(B_{0}\cap B_{1},N_{r}(.))^{\ast }=(B_{0}^{\ast }+B_{1}^{\ast
},\left\vert .\right\vert _{r^{\prime }})$ isometrically.
\end{lemma}

\section{Rotund norm in the complex interpolation spaces}

For any Banach space $X,$ let S$_{X}$ be its unit sphere$.$

\begin{definition}
\label{ut}The norm $\left\Vert .\right\Vert $ on a Banach space $X$ is said
to be rotund if $\bigl\Vert\frac{x+y}{2}\bigr\Vert<1$ whenever $x,y\in S_{X}$
are such that x$\neq y.$
\end{definition}

By \cite[Chap.I,Prop.1.3]{D-G-Z} a norm $\bigl\Vert.\bigr\Vert$ on a Banach
space $X$ is rotund if and only if for every x,y$\in X$ such that $2%
\bigl\Vert x\bigr\Vert^{2}+2\bigl\Vert y\bigr\Vert^{2}-\bigl\Vert x+y%
\bigr\Vert^{2}=0,$ one has $x=y$.

\begin{theorem}
\label{cxw}Let $\overline{B}=(B_{0},B_{1})$ be an interpolation couple and
let $\theta \in (0,1).$ If $B_{0}^{\ast }$ is rotund then $B_{\theta }^{\ast
}$ has the same property$.$
\end{theorem}

By a similar argument analogous to that one of \cite[Lemme 5.2]{Da2} one
shows the following lemma:

\begin{lemma}
\label{ii}Let $X$ be a Banach space. The following assertions are
equivalents:

a) The norm on $X^{\ast }$ is rotund$.$

b) For every closed subspace $Y$ of $X,$ the norm on $Y^{\ast }$ is rotund $%
. $

c) For every closed subspace $Y$ of $X,$ there exists a closed subspace $Z$
of $X$ containing $Y$ such that the norm on $Z^{\ast }$ is rotund.
\end{lemma}

Proof of theorem \ref{cxw}: By a similar argument analogous to that one of 
\cite{Da4} we can assume that $(B_{0},B_{1})$ a regular couple of separable
spaces (by using lemma \ref{ii}).

\bigskip Let $a^{\ast },b^{\ast }\in B_{\theta }^{\ast }=(B_{0}^{\ast
},B_{1}^{\ast })^{\theta };$ \ By lemma \ref{YO} there exist $G,H\in 
\mathcal{F}_{\ast }^{\infty }(B_{0}^{\ast },B^{\ast })$ such that $G(\theta
)=a^{\ast },H(\theta )=b^{\ast }$,%
\begin{eqnarray*}
&&\bigl\Vert a^{\ast }\bigr\Vert_{(B_{0}^{\ast },B_{1}^{\ast })^{\theta
}}^{2} \\
&=&\mathop{\displaystyle \int}\limits_{\mathbb{R}}\bigl\Vert G_{\ast }(i\tau
)\bigr\Vert_{B_{0}^{\ast }}^{2}Q_{0}(\theta ,\tau )d\tau +%
\mathop{\displaystyle \int}\limits_{\mathbb{R}}\bigl\Vert G_{\ast }(1+i\tau )%
\bigr\Vert_{B_{1}^{\ast }}^{2}Q_{1}(\theta ,\tau )d\tau
\end{eqnarray*}%
and%
\begin{eqnarray*}
&&\bigl\Vert b^{\ast }\bigr\Vert_{(B_{0}^{\ast },B_{1}^{\ast })^{\theta
}}^{2} \\
&=&\mathop{\displaystyle \int}\limits_{\mathbb{R}}\bigl\Vert H_{\ast }(i\tau
)\bigr\Vert_{B_{0}^{\ast }}^{2}Q_{0}(\theta ,\tau )d\tau +%
\mathop{\displaystyle \int}\limits_{\mathbb{R}}\bigl\Vert H_{\ast }(1+i\tau )%
\bigr\Vert_{B_{1}^{\ast }}^{2}Q_{1}(\theta ,\tau )d\tau .
\end{eqnarray*}%
Assume now that $2\bigl\Vert a^{\ast }\bigr\Vert_{_{(B_{0}^{\ast
},B_{1}^{\ast })^{\theta }}}^{2}+2\bigl\Vert b^{\ast }\bigr\Vert%
_{_{(B_{0}^{\ast },B_{1}^{\ast })^{\theta }}}^{2}-\bigl\Vert a^{\ast
}+b^{\ast }\bigr\Vert_{_{(B_{0}^{\ast },B_{1}^{\ast })^{\theta }}}^{2}=0.$
By lemma \ref{JM}%
\begin{eqnarray*}
&&\mathop{\displaystyle \int}\limits_{\mathbb{R}}\bigl\Vert G_{\ast }(i\tau
)+H_{\ast }(i\tau )\bigr\Vert_{B_{0}^{\ast }}^{2}Q_{0}(\theta ,\tau )d\tau +%
\mathop{\displaystyle \int}\limits_{\mathbb{R}}\bigl\Vert G_{\ast }(1+i\tau
)+H_{\ast }(1+i\tau )\bigr\Vert_{B_{1}^{\ast }}^{2}Q_{1}(\theta ,\tau )d\tau
\\
&\geq &\bigl\Vert a^{\ast }+b^{\ast }\bigr\Vert_{(B_{0}^{\ast },B_{1}^{\ast
})^{\theta }}^{2},
\end{eqnarray*}

hence\bigskip $\mathop{\displaystyle \int}\limits_{\mathbb{R}}\left[ 2%
\bigl\Vert G_{\ast }(i\tau )\bigr\Vert_{B_{0}^{\ast }}^{2}+2\bigl\Vert %
H_{\ast }(i\tau )\bigr\Vert_{B_{0}^{\ast }}^{2}-\bigl\Vert G_{\ast }(i\tau
)+H_{\ast }(i\tau )\bigr\Vert_{B_{0}^{\ast }}^{2}\right] Q_{0}(\theta ,\tau
)d\tau $

+$\mathop{\displaystyle \int}\limits_{\mathbb{R}}\left[ 2\bigl\Vert G_{\ast
}(1+i\tau )\bigr\Vert_{B_{1}^{\ast }}^{2}+2\bigl\Vert H_{\ast }(1+i\tau )%
\bigr\Vert_{B_{1}^{\ast }}^{2}-\bigl\Vert G_{\ast }(1+i\tau )+H_{\ast
}(1+i\tau )\bigr\Vert_{B_{1}^{\ast }}^{2}\right] Q_{1}(\theta ,\tau )d\tau $

\bigskip $\leq 2\bigl\Vert a^{\ast }\bigr\Vert_{(B_{0}^{\ast },B_{1}^{\ast
})^{\theta }}^{2}+2\bigl\Vert b^{\ast }\bigr\Vert_{(B_{0}^{\ast
},B_{1}^{\ast })^{\theta }}^{2}-2\bigl\Vert a^{\ast }+b^{\ast }\bigr\Vert%
_{(B_{0}^{\ast },B_{1}^{\ast })^{\theta }}^{2}=0.$

This implies that%
\begin{equation*}
\mathop{\displaystyle \int}\limits_{\mathbb{R}}\left[ 2\bigl\Vert G_{\ast
}(i\tau )\bigr\Vert_{B_{0}^{\ast }}^{2}+2\bigl\Vert H_{\ast }(i\tau )%
\bigr\Vert_{B_{0}^{\ast }}^{2}-\bigl\Vert G_{\ast }(i\tau )+H_{\ast }(i\tau )%
\bigr\Vert_{B_{0}^{\ast }}^{2}\right] Q_{0}(\theta ,\tau )d\tau =0,
\end{equation*}

it results that for almost every $\tau \in \mathbb{R}$%
\begin{equation*}
2\bigl\Vert G_{\ast }(i\tau )\bigr\Vert_{B_{0}^{\ast }}^{2}+2\bigl\Vert %
H_{\ast }(i\tau )\bigr\Vert_{B_{0}^{\ast }}^{2}-\bigl\Vert G_{\ast }(i\tau
)+H_{\ast }(i\tau )\bigr\Vert_{B_{0}^{\ast }}^{2}=0.
\end{equation*}

Since $B_{0}^{\ast }$ is rotund , for almost every $\tau \in \mathbb{R}$ $%
G_{\ast }(i\tau )=H_{\ast }(i\tau ).$ Finally by \cite[Lemma 3.10]{Da4} $%
a^{\ast }=b^{\ast }.\blacksquare $

\section{UR and WUR norms in the complex interpolation spaces}

\begin{definition}
\label{iop}(i)The norm $\left\Vert .\right\Vert $ on a Banach space X is
said to be uniformly rotund (UR for short) if $\bigl\Vert x_{n}-y_{n}%
\bigr\Vert\rightarrow _{n\rightarrow +\infty }0$ whenever $x_{n},y_{n}\in
S_{X},$ $n=0,1,...$ are such that lim$_{n\rightarrow +\infty }\bigl\Vert %
x_{n}+y_{n}\bigr\Vert=2$.

(ii) The norm $\bigl\Vert.\bigr\Vert$ on a Banach space X is said to be
weakly uniformly rotund (WUR \ for short) if $x_{n}-y_{n}\rightarrow
_{n\rightarrow +\infty }0$ in the weak topology of X whenever $%
x_{n},y_{n}\in S_{X},$ $n=1,2,...$ are such that lim$_{n\rightarrow +\infty }%
\bigl\Vert x_{n}+y_{n}\bigr\Vert=2$.
\end{definition}

By \cite[Chap.II,Prop.6.2]{D-G-Z} the norm $\bigl\Vert.\bigr\Vert$ on a
Banach space X is WUR, if and only if $x_{n}-y_{n}\rightarrow _{n\rightarrow
+\infty }0$ in the weak topology of $X$ whenever $x_{n},y_{n}\in X,$ $%
n=1,2,...$ are such that $2\bigl\Vert x_{n}\bigr\Vert^{2}+2\bigl\Vert x_{n}%
\bigr\Vert^{2}-\bigl\Vert x_{n}+y_{n}\bigr\Vert^{2}\rightarrow
_{n\rightarrow +\infty }0$ and $\left\{ x_{n}\right\} $ is bounded.

By a similar argument analogous to that one of \cite[Chap.II,Prop.6.2]{D-G-Z}
we show that the norm $\bigl\Vert.\bigr\Vert$ on a Banach space X is
uniformly rotund if and only if $x_{n}-y_{n}\rightarrow _{n\rightarrow
+\infty }0$ in $X$ whenever $x_{n},y_{n}\in X,$ $n=1,2,...$ are such that $2%
\bigl\Vert x_{n}\bigr\Vert^{2}+2\bigl\Vert x_{n}\bigr\Vert^{2}-\bigl\Vert %
x_{n}+y_{n}\bigr\Vert^{2}\rightarrow _{n\rightarrow +\infty }0$ and $\left\{
x_{n}\right\} $ is bounded.

\begin{remark}
\label{xn}Let (x$_{n})_{n\geq 1},(y_{n})_{n\geq 1}$ be two sequences in a
Banach space $X$ such that $(x_{n})_{n\geq 0}$ is bounded in $X$ and $2%
\bigl\Vert x_{n}\bigr\Vert^{2}+2\bigl\Vert y_{n}\bigr\Vert^{2}-\bigl\Vert %
x_{n}+y_{n}\bigr\Vert^{2}\rightarrow _{n\rightarrow +\infty }0,$ then $%
(y_{n})_{\geq 1}$ is bounded in $X.$
\end{remark}

\begin{theorem}
\label{rty}Let ($B_{0},B_{1})$ be an interpolation couple and let $\theta
\in (0,1)$. Assume that the norm on $B_{0}$ is $UR.$ Then the norm on $%
B_{\theta }^{{}}$ is $UR.$
\end{theorem}

Proof. Let $(a_{n})_{n\geq 1},(b_{n})_{n\geq 1}$ be two $C-$bounded
sequences in $B_{\theta }$. For every $n\geq 1,$ there exist $F_{n},H_{n}\in 
\mathcal{F}(\overline{B})$ satisfying $F_{n}(\theta )=a_{n},H_{n}(\theta
)=b_{n}$, $\bigl\Vert a_{n}\bigr\Vert_{B_{\theta }}^{2}\geq \bigl\Vert F_{n}%
\bigr\Vert_{\mathcal{F}(\overline{B})}^{2}-\frac{1}{4n}$ and $\bigl\Vert %
b_{n}\bigr\Vert_{B_{\theta }}^{2}\geq \bigl\Vert H_{n}\bigr\Vert_{\mathcal{F}%
(\overline{B})}^{2}-\frac{1}{4n}.$

For $j\in \left\{ 0,1\right\} $ let

\begin{eqnarray*}
0 &\leq &S_{n}(j,.) \\
&=&2\bigl\Vert F_{n}(j+i.)\bigr\Vert_{B_{j}}^{2}+2\bigl\Vert H_{n}(j+i.)%
\bigr\Vert_{B_{j}}^{2}-\bigl\Vert F_{n}(j+i.)+H_{n}(j+i.)\bigr\Vert%
_{B_{j}}^{2}.
\end{eqnarray*}

\noindent By \cite[Lemma 4.3.2]{Ber-Lof}, $\bigl\Vert F_{n}+H_{n}\bigr\Vert_{%
\mathcal{F}_{\theta }^{2}(\overline{B})}^{{}}\geq \bigl\Vert%
(F_{n}+H_{n})(\theta )\bigr\Vert_{B_{\theta }}^{{}}=\bigl\Vert a_{n}+b_{n}%
\bigr\Vert_{B_{\theta }}^{{}}$. Thus 
\begin{eqnarray*}
&&\mathop{\displaystyle \int}\limits_{\mathbb{R}}S_{n}(0,\tau )Q_{0}(\theta
,\tau )d\tau +\mathop{\displaystyle \int}\limits_{\mathbb{R}}S_{n}(1,\tau
)Q_{1}(\theta ,\tau )d\tau \\
&=&2\bigl\Vert F_{n}\bigr\Vert_{\mathcal{F}_{\theta }^{2}(\overline{B}%
)}^{2}+2\bigl\Vert H_{n}\bigr\Vert_{\mathcal{F}_{\theta }^{2}(\overline{B}%
)}^{2}-\bigl\Vert F_{n}+H_{n}\bigr\Vert_{\mathcal{F}_{\theta }^{2}(\overline{%
B})}^{2} \\
&\leq &2\bigl\Vert F_{n}\bigr\Vert_{\mathcal{F}(\overline{B})}^{2}+2%
\bigl\Vert H_{n}\bigr\Vert_{\mathcal{F}(\overline{B})}^{2}-\bigl\Vert %
F_{n}+H_{n}\bigr\Vert_{\mathcal{F}_{\theta }^{2}(\overline{B})}^{2} \\
&\leq &\frac{1}{n}+2\bigl\Vert a_{n}\bigr\Vert_{B_{\theta }}^{2}+2\bigl\Vert %
b_{n}\bigr\Vert_{B_{\theta }}^{{}}-\bigl\Vert a_{n}+b_{n}\bigr\Vert%
_{B_{\theta }}^{2}.
\end{eqnarray*}

\noindent Assume now that $2$ $\bigl\Vert a_{n}\bigr\Vert_{B_{\theta }}^{2}+2%
\bigl\Vert b_{n}\bigr\Vert_{B_{\theta }}^{2}-\bigl\Vert a_{n}+b_{n}\bigr\Vert%
_{B_{\theta }}^{2}\rightarrow _{n\rightarrow \infty }0.$ Since $%
S_{n}(j,.)\geq 0,$ for every $j\in \left\{ 0,1\right\} $, $%
\mathop{\displaystyle \int}\limits_{\mathbb{R}}S_{n}(0,\tau )Q_{0}(\theta
,\tau )d\tau \rightarrow _{n\rightarrow \infty }0$. We deduce that there is
a subsquence $(S_{n_{k}}(0,.)_{k\geq 0}$ such that $a.s.$ $%
S_{n_{k}}(0,.)\rightarrow _{n\rightarrow \infty }0.$ The norm on $B_{0}$ is $%
UR$, \ it results that a.e $F_{n_{k}}(i.)-H_{n_{k}}(i.)\rightarrow
_{k\rightarrow \infty }0,$ in $B_{0}.$ By \cite[Lemme 4.3.2]{Ber-Lof}%
\begin{eqnarray*}
&&\bigl\Vert a_{n_{k}}-b_{n_{k}}\bigr\Vert_{B_{\theta }}^{{}} \\
&\leq &\left[ \frac{1}{1-\theta }\mathop{\displaystyle \int}\limits_{\mathbb{%
R}}\bigl\Vert F_{n_{k}}(i\tau )-H_{n_{k}}(i\tau )\bigr\Vert%
_{B_{0}}Q_{0}(\theta ,\tau )d\tau \right] ^{\frac{1}{1-\theta }}\times \\
&&\left[ \frac{1}{\theta }\mathop{\displaystyle \int}\limits_{\mathbb{R}}%
\bigl\Vert F_{n_{k}}(1+i\tau )-H_{n_{k}}(1+i\tau )\bigr\Vert%
_{B_{1}}Q_{1}(\theta ,\tau )d\tau \right] ^{\frac{1}{\theta }} \\
&\leq &C^{\prime }\left[ \frac{1}{1-\theta }\mathop{\displaystyle \int}%
\limits_{\mathbb{R}}\bigl\Vert F_{n_{k}}(i\tau )-H_{n_{k}}(i\tau )\bigr\Vert%
_{B_{0}}Q_{0}(\theta ,\tau )d\tau \right] ^{\frac{1}{1-\theta }},
\end{eqnarray*}

where $C^{\prime }$ is a positive constant. By using the Lebesgue dominated
convergence theorem, we obtain that $a_{n_{k}}-b_{n_{k}}\rightarrow
_{k\rightarrow \infty }0$ in $B_{\theta }.$ By a standard argument, we see
that $a_{n}-b_{n}\underset{n\rightarrow \infty }{\rightarrow }0,$ in $%
B_{\theta }$.$\blacksquare $

\section{$LUR$ and weakly $LUR$ norms in complex interpolation spaces}

\begin{definition}
\label{yu}(i)The norm $\bigl\Vert.\bigr\Vert$ on a Banach space X is said to
be locally uniformly rotund (LUR for short) if $\bigl\Vert x_{n}-x\bigr\Vert%
\rightarrow _{n\rightarrow +\infty }0$ whenever $x_{n},x\in X$, $n=1,2,...$
are such that lim$_{n\rightarrow +\infty }\bigl\Vert x_{n}\bigr\Vert=%
\bigl\Vert x\bigr\Vert$ and lim$_{n\rightarrow +\infty }\bigl\Vert x_{n}+x%
\bigr\Vert=2$.

(ii) The norm $\bigl\Vert.\bigr\Vert$ on a Banach space X is said to be
weakly locally uniformly rotund ( weakly LUR \ for short) if $%
x_{n}-x\rightarrow _{n\rightarrow +\infty }0$ in the weak topology of X
whenever $x_{n},x\in X$, $n=1,2,...$ such that lim$_{n\rightarrow +\infty }%
\bigl\Vert x_{n}\bigr\Vert=\bigl\Vert x\bigr\Vert$ and lim$_{n\rightarrow
+\infty }\bigl\Vert x_{n}+x\bigr\Vert=2$.
\end{definition}

\bigskip By \cite[Chap.II,Prop.1.2]{D-G-Z}, the norm $\bigl\Vert.\bigr\Vert$
on a Banach space X is $LUR$ (resp. weakly $LUR)$ if and only if, for all
sequence $(x_{n})_{n\geq 1}$ in $X$ and all $x\in X$ such that 
\begin{equation*}
2\bigl\Vert x_{n}\bigr\Vert^{2}+2\bigl\Vert x\bigr\Vert^{2}-\bigl\Vert %
x_{n}+x\bigr\Vert^{2}\rightarrow _{n\rightarrow +\infty }0,
\end{equation*}

then $\bigl\Vert x_{n}-x\bigr\Vert\rightarrow _{n\rightarrow +\infty }0$
(resp. $x_{n}-x\rightarrow _{n\rightarrow +\infty }0$ in the weak topology
of $X)$.

\begin{theorem}
\label{UT}\cite{Da2}Let $(B_{0},B_{1})$ be an interpolation couple and let $%
\theta \in (0,1).$ Assume that the norm on $B_{0}^{\ast }$ is weakly $LUR.$
Then $B_{\theta }^{\ast }$ has the same property.
\end{theorem}

By a similar argument analogous to that one of \cite[Th.5.3]{Da2} we show:

\begin{theorem}
\label{ro}Let $(B_{0},B_{1})$ be an interpolation couple and let $\theta \in
(0,1).$ Assume that the norm on $B_{0}^{\ast }$ is $LUR.$ Then $B_{\theta
}^{\ast }$ has the same property.
\end{theorem}

\section{Rotund norm in the real interpolation spaces}

The real interpolation spaces have been developed by J.L.Lions and J.Peetre 
\cite{Lio}-\cite{Lio-Peet}. Let $p\in \left[ 1,+\infty \right[ .$ The real
interpolation space $B_{\theta ,p}$ \ is defined by

\begin{equation*}
B_{\theta ,p}=\biggl\{a\in B_{0}+B_{1};\mathrm{\ }\bigl\Vert a\bigr\Vert%
_{B_{\theta ,p}}^{{}}=\left[ \mathop{\displaystyle \int}\limits_{\mathbb{R}%
^{+}}\left[ \frac{K(t,a)}{t^{\theta }}\right] ^{p}\frac{dt}{t}\right] ^{%
\frac{1}{p}}\biggr\}<+\infty
\end{equation*}

\noindent where

\begin{equation*}
K(t,a,\overline{B})=K(t,a)=\inf \left\{ \bigl\Vert a_{0}\bigr\Vert%
_{B_{0}^{{}}}^{{}}+t\bigl\Vert a_{1}\bigr\Vert_{B_{1}^{{}}}^{{}};\mathrm{\ }%
a=a_{0}+a_{1},\mathrm{\ }a_{j}\in B_{j}\mathrm{,}\text{ }j=0,1\right\} .
\end{equation*}

($B_{\theta ,p}$, $\bigl\Vert.\bigr\Vert_{B_{\theta ,p}}^{{}})$ is a Banach
space \cite[Th.3.4.2]{Ber-Lof}$.$

The dual of $B_{\theta ,p}$ will be denoted by $B_{\theta ,p}^{\ast },$ $%
0<\theta <1,$ $1<p<+\infty .$

\noindent The real interpolation space $B_{\theta ,p}$ have the following
properties:

(i) $B_{0}\cap B_{1}$ is dense in $B_{\theta ,p}$ \cite[Th.3.4.2]{Ber-Lof}.

(ii) The dual of $B_{\theta ,p},$ $1<p<\infty ,$ isomorphically identified
with $(B_{0}^{\ast },B_{1}^{\ast })_{\theta ,p^{\prime }}$ if $(B_{0},B_{1})$
is a regular couple \cite[Th.3.7.1]{Ber-Lof}.

(iii) \ For all $p\geq 1$ let 
\begin{equation*}
K_{p}(t,a,\overline{B})=K_{p}(t,a)=\inf \left\{ (\bigl\Vert a_{0}\bigr\Vert%
_{B_{0}^{{}}}^{p}+t^{p}\bigl\Vert a_{1}\bigr\Vert_{B_{0}^{{}}}^{p})^{\frac{1%
}{p}};\mathrm{\ }a=a_{0}+a_{1},\mathrm{\ }a_{j}\in B_{j}\mathrm{,}\text{ }%
j=0,1\right\}
\end{equation*}

and $n_{B_{\theta ,p}}(a)=\left[ \mathop{\displaystyle \int}\limits_{\mathbb{%
R}^{+}}(\frac{K_{p}(t,a)}{t^{\theta }})^{p}\frac{dt}{t}\right] ^{\frac{1}{p}%
}.$ By \cite[Exercice 3.13]{Ber-Lof}, $n_{B_{\theta ,p}}(.)$ is an
equivalent norm on $B_{\theta ,p}.$

\bigskip

\begin{lemma}
\label{nnn}Let $(B_{0},B_{1})$ be a regular couple.

a) If the norms on $B_{0}^{\ast }$ and $B_{1}^{\ast }$ are rotund, then $%
((B_{0}\cap B_{1})^{\ast },N_{2}(.))^{\ast }$ is rotund.

b) If the norm on $B_{0}^{\ast }$ is rotund, then $(B_{0}+B_{1},\left\vert
.\right\vert _{2})^{\ast }$ is rotund.
\end{lemma}

Proof. a) Assume that the norms on $B_{0}^{\ast }$ and $B_{1}^{\ast }$ are
rotund. As $((B_{0}\cap B_{1})^{\ast },N_{2}(.))^{\ast }=(B_{0}^{\ast
}+B_{1}^{\ast },\left\vert .\right\vert _{2})$ by lemma \ref{Dual}, it
suffices to show that $(B_{0}^{\ast }+B_{1}^{\ast },\left\vert .\right\vert
_{2})$ is rotund. For that, let $a^{\ast },b^{\ast }\in B_{0}^{\ast
}+B_{1}^{\ast }$ such that $2\left\vert a^{\ast }\right\vert
_{2}^{2}+2\left\vert b^{\ast }\right\vert _{2}^{2}-\left\vert a^{\ast
}+b^{\ast }\right\vert _{2}^{2}=0.$ There exist $x_{j}^{\ast },y_{j}^{\ast
}\in B_{j}^{\ast }$ and $j\in \left\{ 0,1\right\} $ satisfying $a^{\ast
}=x_{0}^{\ast }+x_{1}^{\ast },$ $b^{\ast }=y_{0}^{\ast }+y_{1}^{\ast }$, $%
\left\vert a^{\ast }\right\vert _{2}^{2}=\bigl\Vert x_{0}^{\ast }\bigr\Vert%
_{B_{0}^{\ast }}^{2}+\bigl\Vert x_{1}^{\ast }\bigr\Vert_{B_{1}^{\ast }}^{2}$
and

$\left\vert b^{\ast }\right\vert _{2}^{2}=\bigl\Vert y_{0}^{\ast }\bigr\Vert%
_{B_{0}^{\ast }}^{2}+\bigl\Vert y_{1}^{\ast }\bigr\Vert_{B_{1}^{\ast }}^{2}.$
It results that $2\bigl\Vert x_{0}^{\ast }\bigr\Vert_{B_{0}^{\ast }}^{2}+2%
\bigl\Vert y_{0}^{\ast }\bigr\Vert_{B_{0}^{\ast }}^{2}-\bigl\Vert %
x_{0}^{\ast }+y_{0}^{\ast }\bigr\Vert_{B_{0}^{\ast }}^{2}+2\bigl\Vert %
x_{1}^{\ast }\bigr\Vert_{B_{1}^{\ast }}^{2}+2\bigl\Vert y_{1}^{\ast }%
\bigr\Vert_{B_{1}^{\ast }}^{2}-\bigl\Vert x_{1}^{\ast }+y_{1}^{\ast }%
\bigr\Vert_{B_{1}^{\ast }}^{2}\leq 2\left\vert a^{\ast }\right\vert
_{2}^{2}+2\left\vert b^{\ast }\right\vert _{2}^{2}-\left\vert a^{\ast
}+b^{\ast }\right\vert _{2}^{2}=0.$ This implies that $2\bigl\Vert %
x_{0}^{\ast }\bigr\Vert_{B_{0}^{\ast }}^{2}+2\bigl\Vert y_{0}^{\ast }%
\bigr\Vert_{B_{0}^{\ast }}^{2}-\bigl\Vert x_{0}^{\ast }+y_{0}^{\ast }%
\bigr\Vert_{B_{0}^{\ast }}^{2}=0$ and $2\bigl\Vert x_{1}^{\ast }\bigr\Vert%
_{B_{1}^{\ast }}^{2}+2\bigl\Vert y_{1}^{\ast }\bigr\Vert_{B_{1}^{\ast }}^{2}-%
\bigl\Vert x_{1}^{\ast }+y_{1}^{\ast }\bigr\Vert_{B_{1}^{\ast }}^{2}=0.$
Since the norms on $B_{0}^{\ast },B_{1}^{\ast }$ are rotund, $x_{j}^{\ast
}=y_{j}^{\ast },$ $j\in \left\{ 0,1\right\} ,$ i.e. $a^{\ast }=b^{\ast
}.\blacksquare $

b) Assume that the norm on $B_{0}^{\ast }$ is rotund. By lemma \ref{Dual},
it is enough to show that $(B_{0}^{\ast }\cap B_{1}^{\ast },N_{2}(.))$ is
rotund. Let $a^{\ast },b^{\ast }\in B_{0}^{\ast }\cap B_{1}^{\ast }$ such
that $2N_{2}(a^{\ast })^{2}+2N_{2}(b^{\ast })^{2}-N_{2}(a^{\ast }+b)^{2}=0.$
Observe that $2\bigl\Vert a_{{}}^{\ast }\bigr\Vert_{B_{0}^{\ast }}^{2}+2%
\bigl\Vert b_{{}}^{\ast }\bigr\Vert_{B_{0}^{\ast }}^{2}-\bigl\Vert %
a_{{}}^{\ast }+b^{\ast }\bigr\Vert_{B_{0}^{\ast }}^{2}+2\bigl\Vert %
a_{{}}^{\ast }\bigr\Vert_{B_{1}^{\ast }}^{2}+2\bigl\Vert b_{{}}^{\ast }%
\bigr\Vert_{B_{1}^{\ast }}^{2}-\bigl\Vert a_{{}}^{\ast }+b^{\ast }\bigr\Vert%
_{B_{1}^{\ast }}^{2}\leq 2N_{2}(a^{\ast })^{2}+2N_{2}(b^{\ast
})^{2}-N_{2}(a^{\ast }+b^{\ast })^{2}=0.$ It follows that $2\bigl\Vert %
a_{{}}^{\ast }\bigr\Vert_{B_{0}^{\ast }}^{2}+2\bigl\Vert b_{{}}^{\ast }%
\bigr\Vert_{B_{0}^{\ast }}^{2}-\bigl\Vert a_{{}}^{\ast }+b^{\ast }\bigr\Vert%
_{B_{0}^{\ast }}^{2}=0$. Since $B_{0}^{\ast }$ is rotund, $a^{\ast }=b^{\ast
}.\blacksquare $

\begin{corollary}
\bigskip \label{cc}Let ($B_{0},B_{1})$ be an interpolation couple, let $%
\theta \in (0,1)$ and let $p\in (1,+\infty )$ . If the norm on $B_{0}^{\ast
} $ is rotund, then $B_{\theta ,p}^{\ast }$ admits an equivalent rotund dual
norm.
\end{corollary}

Proof. Let $0<\theta _{0}<\theta _{1}<1$ and let $\eta \in (0,1)$ such that $%
\theta =(1-\eta )\theta _{0}+\eta \theta _{1}.$ By the reiteration theorem (%
\cite[Th.4.7.2]{Ber}) $(B_{\theta _{0}},B_{\theta _{1}})_{\eta ,p}=B_{\theta
,p}.$ Now consider the canonical injection $i:B_{\theta _{0}}\cap B_{\theta
_{1}}\rightarrow B_{\theta ,p}$, by recall i) $i(B_{\theta _{0}}\cap
B_{\theta _{1}})$ is dense in $B_{\theta ,p}$, hence $i^{\ast }:B_{\theta
,p}^{\ast }\rightarrow (B_{\theta _{0}}\cap B_{\theta _{1}})^{\ast }$ is
continuous and one to one. On the other hand by using theorem \ref{cxw} and
lemma \ref{nnn}, $\left[ (B_{\theta _{0}}\cap B_{\theta _{1}})^{{}},N_{2}(.)%
\right] ^{\ast }$ is rotund (note that by recall 2 bis, $(B_{\theta
_{0}},B_{\theta _{1}})$ is a regular couple), hence by using \cite[%
Chap.II,Th.2.4]{D-G-Z}, we obtain that $B_{\theta ,p}^{\ast }$ admits an
equivalent rotund dual norm.$\blacksquare $

\begin{remark}
\label{Ug}In \cite{Da3} we show that, there exists an interpolation couple $%
(B_{0},B_{1})$ such that $B_{0}^{{}}$ is a subspace of $B_{1},$ $B_{0}^{\ast
}$ admits an equivalent $LUR$ norm (hence this norm is rotund) but $%
B_{\theta }^{\ast },B_{\theta ,p}^{\ast }$ not admit any equivalent rotund
norms. Thus the condition that the dual norm of $B_{0}^{\ast }$ is rotund is
necessary.
\end{remark}

\section{UR and WUR norms in the real interpoltion spaces}

\begin{theorem}
\label{Day}Let $(B_{0},B_{1})$ be an interpolation couple, let $0<\theta <1$
and let $1<p<+\infty .$ Assume that the norm on $B_{0}$ is $UR,$ then $%
B_{\theta ,p}$ admits an equivalent $UR$ norm.
\end{theorem}

\begin{remark}
\label{P}On $\mathbb{R}^{+}$ denotes $d\mu =\frac{dt}{t^{2\theta +1}}.$ As $%
\mu $ is $\sigma -finite$ measure$,$ there exists a sequence $(\Omega
_{n})_{n\geq 1}$ of pairwise disjoint measurable subsets of $\ \mathbb{R}%
^{+} $ such that $\mathbb{R}^{+}=\underset{n\geq 1}{%
\mathop{\displaystyle
\sum }}\Omega _{n}$ and $\mu (\Omega _{n})<+\infty ,$ for all $n\geq 1.$ One
defines the measure $\mu _{1}$ by $\mu _{1}$($E)=\underset{n\geq 1}{%
\mathop{\displaystyle \sum }}\frac{\mu (E\cap \Omega _{n})}{\mu (\Omega
_{n})2^{n}},$ where $E$ is a measurable subset of $\ \mathbb{R}^{+}.$ It is
clear that $\mu _{1}$ is a probability measure. One defines the operator $%
U:L^{2}(\mu )\rightarrow L^{2}(\mu _{1}),$ by $U(f)(\omega )=2^{n}\mu
(\Omega _{n})f(\omega ),$ if $\omega \in \Omega _{n}.$ we remark that $U$ is
a an isometric operator and $U\otimes I_{X}:L^{2}(\mu ,X)\rightarrow
L^{2}(\mu _{1},X)$ is again an isometric operator for any Banach $X.$
\end{remark}

\begin{lemma}
\label{ci}Let $(B_{0},B_{1})$ be an interpolation couple and let $0<\theta
<1.$ Assume that the norm on $B_{0}$ is $UR,$ then $B_{\theta ,2}$ admits an
equivalent $UR$ norm$.$
\end{lemma}

Proof. By the reiteration theorem \cite[Th.4.7.2]{Ber-Lof}, $W=(B_{\theta
_{0}},B_{\theta _{1}})_{\eta ,2}=B_{\theta ,2}$ isomorphically$,$ where $%
0<\theta _{0}<\theta _{1}<1,$ $0<\eta <1$ and $\theta =(1-\eta )\theta
_{0}+\eta \theta _{1}.$ Let ($a_{n})_{n\geq 1},$ $(b_{n})_{n\geq 1}$ be two
bounded sequences in $W$ such that $2\bigl\Vert a_{n}\bigr\Vert_{W}^{2}+2%
\bigl\Vert b_{n}\bigr\Vert_{W}^{2}-\bigl\Vert a_{n}+b_{n}\bigr\Vert_{W}^{2}%
\underset{n\rightarrow +\infty }{\rightarrow }0,$ where $\bigl\Vert.%
\bigr\Vert_{W}^{{}}=n_{X_{\eta ,2}}(.)$ and $X_{j}=B_{\theta _{j}},$ $j\in
\left\{ 0,1\right\} .$ There is a closed separable subspace $Y_{j}$ of $%
X_{j} $ such that the closed subspace spanned by $\left\{ a_{n};\text{ }%
n\geq 1\right\} \cup $ $\left\{ b_{n};\text{ }n\geq 1\right\} $ continuously
embeds in $Z=(Y_{0},Y_{1})_{\eta ,2},$ $j\in \left\{ 0,1\right\} $, $%
\left\Vert a_{n}\right\Vert _{W}=\left\Vert a_{n}\right\Vert _{Z}$ and $%
\bigl\Vert b_{n}\bigr\Vert_{W}^{{}}=\bigl\Vert b_{n}\bigr\Vert_{Z}^{{}},$ $%
n\geq 1.$ We deduce that $2\bigl\Vert a_{n}\bigr\Vert_{Z}^{2}+2\bigl\Vert %
b_{n}\bigr\Vert_{Z}^{2}-\bigl\Vert a_{n}+b_{n}\bigr\Vert_{Z}^{2}\underset{%
n\rightarrow +\infty }{\rightarrow }0,$ here

$\bigl\Vert.\bigr\Vert_{Z}^{{}}=n_{Y_{\eta ,2}}(.).$

Since $Y_{0},Y_{1}$ are polish spaces, by selection theorem, there exist two
sequences $(a_{n}^{j})_{n\geq 1},(b_{n}^{j})_{n\geq 1}$ in $L^{2}(Y_{j}),$ $%
j\in \left\{ 0,1\right\} $ such that $%
(K_{2}^{{}}(t,a_{n}))^{2}=(K_{2}^{{}}(t,a_{n},Y_{0},Y_{1}))^{2}$ \TEXTsymbol{%
>}$\bigl\Vert a_{n}^{0}(t)\bigr\Vert_{Y_{0}}^{2}+t^{2}\bigl\Vert a_{n}^{1}(t)%
\bigr\Vert_{Y_{1}}^{2}-\frac{h(t)}{n},$ $a_{n}=a_{n}^{0}(t)+a_{n}^{1}(t)$, ($%
K_{2}^{{}}(t,b_{n}))^{2}=(K_{2}^{{}}(t,b_{n},Y_{0},Y_{1}))^{2}$ \TEXTsymbol{>%
}$\left\Vert b_{n}^{0}(t)\right\Vert _{Y_{0}}^{2}+t^{2}\left\Vert
b_{n}^{1}(t)\right\Vert _{Y_{1}}^{2}-\frac{h(t)}{n}$ and $%
b_{n}=b_{n}^{0}(t)+b_{n}^{1}(t)$, for almost every $t\in \mathbb{R}^{+}$ and
all $n\geq 1,$ where $\ h$ $\in L^{1}(\mu )$ and $h(t)>0,$ for almost every $%
t\in \mathbb{R}^{+}.$ We observe that $2\bigl\Vert a_{n}^{0}\bigr\Vert%
_{L^{2}(\mu ,Y_{0})}^{2}+2\bigl\Vert b_{n}^{0}\bigr\Vert_{L^{2}(\mu
,Y_{0})}^{2}-\bigl\Vert a_{n}^{0}+b_{n}^{0}\bigr\Vert_{L^{2}(\mu ,Y_{0})}^{2}%
\underset{n\rightarrow +\infty }{\rightarrow }0$ and $2\bigl\Vert ta_{n}^{1}%
\bigr\Vert_{L^{2}(\mu ,Y_{1})}^{2}+2\bigl\Vert tb_{n}^{1}\bigr\Vert%
_{L^{2}(\mu ,Y_{1})}^{2}-\bigl\Vert ta_{n}^{1}+tb_{n}^{1}\bigr\Vert%
_{L^{2}(\mu ,Y_{1})}^{2}\underset{n\rightarrow +\infty }{\rightarrow }0.$ By
theorem \ref{rty}, $X_{0},X_{1}$ are $UR,$ hence $Y_{0},Y_{1}$ are $UR.$ On
the other hand by remark \ref{P} and the result of \cite{Day} the norm on $%
L^{2}(\mu ,Y_{j})$ is $UR,$ $j\in \left\{ 0,1\right\} ,$ this implies that $%
\bigl\Vert a_{n}^{0}-b_{n}^{0}\bigr\Vert_{L^{2}(\mu ,Y_{0})}^{2}\underset{}{%
\rightarrow _{n\rightarrow +\infty }}0$ and $\bigl\Vert ta_{n}^{1}-tb_{n}^{1}%
\bigr\Vert_{L^{2}(\mu ,Y_{1})}^{2}\underset{n\rightarrow +\infty }{%
\rightarrow }0,$ hence $\bigl\Vert a_{n}-b_{n}\bigr\Vert_{Z}^{2}\rightarrow
_{n\rightarrow +\infty }0.$ Thus $W$ is $UR.\blacksquare $

Proof of theorem \ref{Day}. By using the reiteration theorem, we obtain that$%
(B_{\theta _{0},2},B_{\theta _{1},q})_{\eta }=B_{\theta ,p}$ isomorphically$%
, $ where $0<\theta _{0}<\theta _{1}<1,$ $0<\eta <1,$ $\theta =(1-\eta
)\theta _{0}+\eta \theta _{1}$ and $\frac{1}{p}=\frac{1-\eta }{2}+\frac{\eta 
}{q}.$ On the other hand by lemma \ref{ci}, $B_{\theta _{0},2}$ is $UR$ and
by theorem \ref{rty} $(B_{\theta _{0},2},B_{\theta _{1},q})_{\eta
}=B_{\theta ,p}$ is $UR.\blacksquare $

\begin{theorem}
\label{xxx}Let ($B_{0},B_{1})$ be an interpolation couple, let $\theta \in
(0,1)$ and let $p\in (1,+\infty )$. Assume that the norm on $B_{0}^{{}}$ is $%
WUR,$ then $B_{\theta ,p}^{{}}$ admis an equivalent $WUR$ norm.
\end{theorem}

\begin{lemma}
\label{bv}Let $(B_{0},B_{1})$ be an interpolation couple. Assume that the $%
B_{0}^{{}}$,$B_{1}^{{}}$ are $WUR$. Then $(B_{0}^{{}}+B_{1}^{{}},\left\vert
.\right\vert _{2})$ is $WUR$.
\end{lemma}

Proof. Let ($a_{n})_{n\geq 1},$ $(b_{n})_{n\geq 1}$ be two sequences in $%
B_{0}+B_{1}$ such that ($a_{n})_{n\geq 1}$ is bounded in $B_{0}+B_{1}$ and $%
2\left\vert a_{n}\right\vert _{2}^{2}+2\left\vert b_{n}\right\vert
_{2}^{2}-\left\vert a_{n}+b_{n}\right\vert _{2}^{2}\underset{n\rightarrow
+\infty }{\rightarrow }0.$ For every $n\geq 1,$ there exist $%
x_{n}^{j},y_{n}^{j},\in B_{j},$ $j\in \left\{ 0,1\right\} $ such that $%
a_{n}=x_{n}^{0}+x_{n}^{1},$ $b_{n}=y_{n}^{0}+y_{n}^{1},$ $\left\vert
a_{n}\right\vert ^{2}>\bigl\Vert x_{n}^{0}\bigr\Vert_{B_{0}}^{2}+\bigl\Vert %
x_{n}^{1}\bigr\Vert_{B_{1}}^{2}-\frac{1}{4n}$ and $\left\vert
b_{n}\right\vert ^{2}>\bigl\Vert y_{n}^{0}\bigr\Vert_{B_{0}}^{2}+\bigl\Vert %
y_{n}^{1}\bigr\Vert_{B_{1}}^{2}-\frac{1}{4n}.$ Note that%
\begin{eqnarray*}
2\left\vert a_{n}\right\vert _{2}^{2}+2\left\vert b_{n}\right\vert
_{2}^{2}-\left\vert a_{n}+b_{n}\right\vert _{2}^{2} &\geq & \\
&&2\bigl\Vert x_{n}^{0}\bigr\Vert_{B_{0}}^{2}+2\bigl\Vert y_{n}^{0}\bigr\Vert%
_{B_{0}}^{2}-\bigl\Vert x_{n}^{0}+y_{n}^{0}\bigr\Vert_{B_{0}}^{2}+2%
\bigl\Vert x_{n}^{1}\bigr\Vert_{B_{1}}^{2} \\
&&+2\bigl\Vert y_{n}^{1}\bigr\Vert_{B_{1}}^{2}-\bigl\Vert x_{n}^{1}+y_{n}^{1}%
\bigr\Vert_{B_{1}}^{2}-\frac{1}{n}.
\end{eqnarray*}

It follows that $2\bigl\Vert x_{n}^{0}\bigr\Vert_{B_{0}}^{2}+2\bigl\Vert %
y_{n}^{0}\bigr\Vert_{B_{0}}^{2}-\bigl\Vert x_{n}^{0}+y_{n}^{0}\bigr\Vert%
_{B_{0}}^{2}\underset{n\rightarrow +\infty }{\rightarrow }0$ and $2%
\bigl\Vert
x_{n}^{1}\bigr\Vert_{B_{1}}^{2}+2\bigl\Vert y_{n}^{1}\bigr\Vert_{B_{1}}^{2}-%
\bigl\Vert x_{n}^{1}+y_{n}^{1}\bigr\Vert_{B_{1}}^{2}\underset{n\rightarrow
+\infty }{\rightarrow }0.$ Since $B_{0},B_{1}$ are $WUR,$ $%
x_{n}^{j}-y_{n}^{j}\rightarrow 0$ in the weak topology of $B_{j},$ $j\in
\left\{ 0,1\right\} .$ It results that $%
a_{n}-b_{n}=x_{n}^{0}-y_{n}^{0}+x_{n}^{1}-y_{n}^{1}\underset{n\rightarrow
+\infty }{\rightarrow }0$ in the weak topology of $B_{0}+B_{1}.\blacksquare $

Proof of theorem \ref{xxx}. We can assume that $(B_{0},B_{1})$ is a regular
couple. Let $0<\theta _{0}<\theta _{1}<1$ and let $\eta \in (0,1)$ such that 
$\theta =(1-\eta )\theta _{0}+\eta \theta _{1};$ by the reiteration theorem (%
\cite[Th.4.7.2]{Ber-Lof}) $(B_{\theta _{0}},B_{\theta _{1}})_{\eta
,p}=B_{\theta ,p}.$ Consider now the canonical injection $i:B_{\theta
,p}\rightarrow B_{\theta _{0}}+B_{\theta _{1}}.$ By \cite[Th.4.1]{Da2} $%
B_{\theta _{0}},B_{\theta _{1}}$ are $WUR,$ by using lemma \ref{bv} we
obtain that ($B_{\theta _{0}}+B_{\theta _{1}},\left\vert .\right\vert _{2})$
is $WUR.$ On the other hand by \cite[ Th.2.7.1]{Ber-Lof}, $(B_{\theta
_{0}}+B_{\theta _{1}})^{\ast }=B_{\theta _{0}}^{\ast }\cap B_{\theta
_{1}}^{\ast }.$ But by recall (ii), $(B_{\theta _{0}}^{\ast },B_{\theta
_{1}}^{\ast })_{\eta ,p^{\prime }}=B_{\theta ,p}^{\ast },$ this implies that 
$i^{\ast }(B_{\theta _{0}}^{\ast }\cap B_{\theta _{1}}^{\ast })$ is dense in 
$B_{\theta ,p}^{\ast }$ (by recall (i)) ($i^{\ast }:$($B_{\theta
_{0}}+B_{\theta _{1}}$)$^{\ast }=B_{\theta _{0}}^{\ast }\cap B_{\theta
_{1}}^{\ast }\rightarrow B_{\theta ,p}^{\ast }$)$.$ Thus by \cite[%
Chap.II,Th.6.8]{D-G-Z} $B_{\theta ,p}^{{}}$ admits an equivalent $WUR$ norm$%
.\blacksquare $

\section{norms $LUR$ and weakly $LUR$ in the real interpolated spaces}

\begin{theorem}
\label{or}Let $(B_{0},B_{1})$ be an interpolation space, let $\theta \in
(0,1)$ and let $p\in (1,+\infty ).$ Assume that $B_{0}^{\ast }$ is $LUR.$
Then $B_{\theta ,p}^{\ast }$ admits an equivalent $LUR$ dual norm$.$
\end{theorem}

\begin{lemma}
\label{bk}Let $(B_{0},B_{1})$ be a regular couple. Assume that $B_{0}^{\ast
} $ and $B_{1}^{\ast }$ are $LUR$. Then $(B_{0}+B_{1},\left\vert
.\right\vert _{2})^{\ast }$ is $LUR.$
\end{lemma}

Proof. By lemma \ref{Dual}, it suffices to show that $(B_{0}^{\ast }\cap
B_{1}^{\ast },.N_{2}(.))$ is $LUR.$ Let $(a_{n}^{\ast })_{n\geq 1}$ be a
sequence in $(B_{0}^{\ast }\cap B_{1}^{\ast },.N_{2}(.))$ and let $a^{\ast
}\in B_{0}^{\ast }\cap B_{1}^{\ast }$ such that $2N_{2}(a_{n}^{\ast
})^{2}+2N_{2}(a^{\ast })^{2}-N_{2}(a_{n}^{\ast }+a^{\ast })^{2}=2\bigl\Vert %
a_{n}^{\ast }\bigr\Vert_{B_{0}}^{2}+2\bigl\Vert a_{{}}^{\ast }\bigr\Vert%
_{B_{0}}^{2}-\bigl\Vert a_{n}^{\ast }+a^{\ast }\bigr\Vert_{B_{0}}^{2}+2%
\bigl\Vert a_{n}^{\ast }\bigr\Vert_{B_{1}}^{2}+2\bigl\Vert a_{{}}^{\ast }%
\bigr\Vert_{B_{1}}^{2}-\bigl\Vert a_{n}^{\ast }+a_{{}}^{\ast }\bigr\Vert%
_{B_{1}}^{2}\underset{n\rightarrow +\infty }{\rightarrow }0,$ it results
that $2\bigl\Vert a_{n}^{\ast }\bigr\Vert_{B_{j}}^{2}+2\bigl\Vert %
a_{{}}^{\ast }\bigr\Vert_{B_{j}}^{2}-\bigl\Vert a_{n}^{\ast }+a^{\ast }%
\bigr\Vert_{B_{j}}^{2}\underset{n\rightarrow +\infty }{\rightarrow }0,$
since $B_{j}^{\ast }$ is $LUR,$ $a_{n}^{\ast }\rightarrow _{n\rightarrow
+\infty }a^{\ast }$ in $B_{j}^{\ast },$ $j\in \left\{ 0,1\right\} .$ One
deduces that $a_{n}^{\ast }\rightarrow _{n\rightarrow +\infty }a^{\ast }$ in 
$B_{0}^{\ast }\cap B_{1}^{\ast }.\blacksquare $

Proof of theorem \ref{or}. We can assume that $(B_{0},B_{1})$ is a regular
couple. Let $\theta _{0}<\theta _{1}<1$ and let $0<\eta <1$ such that $%
\theta =(1-\eta )\theta _{0}+\eta \theta _{1}.$ By the reiteration theorem (%
\cite[Th.4.7.2]{Ber-Lof}) $(B_{\theta _{0}},B_{\theta _{1}})_{\eta
,p}=B_{\theta ,p}$ $.$ Let $i:B_{\theta ,p}\rightarrow B_{\theta
_{0}}+B_{\theta _{1}}$ be the canonical injection. Throw the proof of
theorem \ref{xxx}, we showed that $i^{\ast }((B_{\theta _{0}}+B_{\theta
_{1}})^{\ast })$ is dense in $B_{\theta ,p}^{\ast }$ ($i^{\ast }:(B_{\theta
_{0}}+B_{\theta _{1}})^{\ast }\rightarrow B_{\theta ,p}^{\ast })$. By using
lemma \ref{bk} et \cite[Chap.II,Th.2.1]{D-G-Z}, we obtain that $B_{\theta
,p}^{\ast }$ admis an equivalent $LUR$ dual norm (note that by recall 2 bis
that ($B_{\theta _{0}},B_{\theta _{1}})$ is regular).$\blacksquare $

\begin{lemma}
\label{gu}Let $(B_{0},B_{1})$ be a regular couple. Assume that $B_{0}^{\ast
},B_{1}^{\ast }$ are weakly $LUR.$ Then ($B_{0}\cap B_{1}$,$N_{2}(.))^{\ast
} $ is weakly $LUR$
\end{lemma}

Proof: By lemma \ref{Dual}, it is enough to show that ($B_{0}^{\ast
}+B_{1}^{\ast },\left\vert .\right\vert _{2})$ is weakly $LUR.$ For that,
let $(a_{n}^{\ast })_{n\geq 1}$ be a sequence in $B_{0}^{\ast }+B_{1}^{\ast
} $ and $a^{\ast }\in B_{0}^{\ast }+B_{1}^{\ast }$ such that $2\left\vert
a_{n}^{\ast }\right\vert _{2}^{2}+2\left\vert a^{\ast }\right\vert
_{2}^{2}-\left\vert a_{n}+a^{\ast }\right\vert _{2}^{2}\rightarrow
_{n\rightarrow +\infty }0.$ For every $n\geq 1,$ there exist $%
x_{n}^{j},c_{j}\in B_{j}^{\ast },$ $j\in \left\{ 0,1\right\} $ such that $%
a_{n}^{\ast }=x_{n}^{0}+x_{n}^{1},$ $a^{\ast }=c_{0}+c_{1},$ $\left\vert
a_{n}^{\ast }\right\vert _{2}^{2}=\bigl\Vert x_{n}^{0}\bigr\Vert%
_{B_{0}^{\ast }}^{2}+\bigl\Vert x_{n}^{1}\bigr\Vert_{B_{1}^{\ast }}^{2},$ $%
\left\vert a_{{}}^{\ast }\right\vert ^{2}=\bigl\Vert c_{0}\bigr\Vert%
_{B_{0}^{\ast }}^{2}+\bigl\Vert c_{1}\bigr\Vert_{B_{1}^{\ast }}^{2}$ $.$
This implies that $2\bigl\Vert x_{n}^{0}\bigr\Vert_{B_{0}^{\ast }}^{2}+2%
\bigl\Vert c_{0}\bigr\Vert_{B_{0}^{\ast }}^{2}-\bigl\Vert x_{n}^{0}+c_{0}%
\bigr\Vert_{B_{0}^{\ast }}^{2}\rightarrow _{n\rightarrow +\infty }0$ and $2%
\bigl\Vert x_{n}^{1}\bigr\Vert_{B_{1}^{\ast }}^{2}+2\bigl\Vert c_{1}%
\bigr\Vert_{B_{1}^{\ast }}^{2}-\bigl\Vert x_{n}^{1}+c_{1}\bigr\Vert%
_{B_{1}^{\ast }}^{2}\rightarrow _{n\rightarrow +\infty }0.$ Since $%
B_{0}^{\ast },B_{1}^{\ast }$ are weakly $LUR,$ $x_{n}^{j}-c_{j}\rightarrow
_{n\rightarrow +\infty }0$ in the weak topology of $B_{j}^{\ast },$ $j\in
\left\{ 0,1\right\} $. Thus $x_{n}-y_{n}\rightarrow _{n\rightarrow +\infty
}0 $ in the weak topology of $B_{0}^{\ast }+B_{1}^{\ast }.\blacksquare $

\begin{theorem}
\label{sz}Let $(B_{0},B_{1})$ be an interpolation couple let $\theta \in
(0,1)$ and let $p\in (1,+\infty )$. Assume that $B_{0}^{\ast }$ is weakly $%
LUR.$ Then $B_{\theta ,p}^{\ast }$ admis an equivalent $LUR$ dual norm$.$
\end{theorem}

Proof. We can assume that $(B_{0},B_{1})$ is a regular couple. Let $\alpha
<\beta <1$ and let $0<\eta <1$ such that $\theta =(1-\eta )\alpha +\eta
\beta .$ By the reiteration theorem $(B_{\alpha },B_{\beta })_{\eta
,p}=B_{\theta ,p}.$ Let $i:B_{\alpha }\cap B_{\beta }\rightarrow B_{\theta
,p}$ be the canonical injection (note that $i^{\ast }:B_{\theta ,p}^{\ast
}\rightarrow (B_{\alpha }\cap B_{\beta })^{\ast }).$ For every $x^{\ast }\in
B_{\theta ,p}^{\ast },$ let $N(x^{\ast })=\left[ \bigl\Vert x^{\ast }%
\bigr\Vert_{B_{\theta ,p}^{\ast }}^{2}+\left\vert i^{\ast }x^{\ast
}\right\vert ^{2}\right] ^{1/2}.$ It suffices to show that ($B_{\theta
,p}^{\ast },N(.))$ is weakly $LUR.$ For that, let ($x_{n}^{\ast })_{n\geq 1}$
be a sequence in $B_{\theta ,p}^{\ast }$ and let $x^{\ast }\in B_{\theta
,p}^{\ast }$ such that $2N(x_{n}^{\ast })^{2}+2N(x^{\ast
})^{2}-N(x_{n}^{\ast }+x^{\ast })^{2}\underset{n\rightarrow +\infty }{%
\rightarrow }0.$ By \cite{Da2} $B_{\alpha }^{\ast },B_{\beta }^{\ast }$ are
weakly $LUR,$ hence ($B_{\alpha }\cap B_{\beta })^{\ast }$ is weakly $LUR$
by lemma \ref{gu}, hence $i^{\ast }x_{n}^{\ast }\rightarrow _{n\rightarrow
+\infty }i^{\ast }x^{\ast }$ in the weak topology of ($B_{\alpha }\cap
B_{\beta })^{\ast }.$ On the other hand by hypothesis $B_{0}^{\ast }$ is
weakly $LUR,$ by \cite{Su}, $B_{0}^{\ast }$ has the Radon-Nikodym property,
by using recall 1) and \cite[Corol.4.5.2]{Ber-Lof}, one obtains that $%
B_{\alpha }^{\ast }=(B_{0}^{\ast },B_{1}^{\ast })^{\alpha }=(B_{0}^{\ast
},B_{1}^{\ast })_{\alpha }$ and $B_{\beta }^{\ast }=(B_{0}^{\ast
},B_{1}^{\ast })^{\beta }=(B_{0}^{\ast },B_{1}^{\ast })_{\beta }.$ Thus by
recall (ii), $B_{\theta ,p}^{\ast }=((B_{0}^{\ast },B_{1}^{\ast })_{\alpha
},(B_{0}^{\ast },B_{1}^{\ast })_{\beta })_{\eta ,p^{\prime }}.$ The recall 2
bis shows us that $(B_{0}^{\ast },B_{1}^{\ast })_{\alpha }\cap (B_{0}^{\ast
},B_{1}^{\ast })_{\beta }$ is dense in $(B_{0}^{\ast },B_{1}^{\ast
})_{\alpha }$ and in $(B_{0}^{\ast },B_{1}^{\ast })_{\beta }$ by recall
(ii), one has that $B_{\theta ,p}^{\ast \ast }=((B_{0}^{\ast },B_{1}^{\ast
})_{\alpha }^{\ast },(B_{0}^{\ast },B_{1}^{\ast })_{\beta }^{\ast })_{\eta
,p}$. It results by recall (i) that $i^{\ast \ast }\left[ (B_{0}^{\ast
},B_{1}^{\ast })_{\alpha }^{\ast }\cap (B_{0}^{\ast },B_{1}^{\ast })_{\beta
}^{\ast }\right] $ is dense in $B_{\theta ,p}^{\ast \ast }.$ On the other
hand by \cite[Th.2.7.1]{Ber-Lof}, $(B_{0}^{\ast },B_{1}^{\ast })_{\alpha
}^{\ast }\cap (B_{0}^{\ast },B_{1}^{\ast })_{\beta }^{\ast }=\left[
(B_{0}^{\ast },B_{1}^{\ast })_{\alpha }^{{}}+(B_{0}^{\ast },B_{1}^{\ast
})_{\beta }^{{}}\right] ^{\ast }=\left[ B_{\alpha }^{\ast }+B_{\beta }^{\ast
}\right] ^{\ast }.$ Note that if $a^{\ast \ast }\in \left[ (B_{0}^{\ast
},B_{1}^{\ast })_{\alpha }^{\ast }\cap (B_{0}^{\ast },B_{1}^{\ast })_{\beta
}^{\ast }\right] ,$ then $\left\langle x_{n}^{\ast },i^{\ast \ast }(a^{\ast
\ast })\right\rangle =\left\langle i^{\ast }(x_{n}^{\ast }),a^{\ast \ast
}\right\rangle \rightarrow _{n\rightarrow +\infty }\left\langle i^{\ast
}(x_{{}}^{\ast }),a^{\ast \ast }\right\rangle =\left\langle x^{\ast
},i^{\ast \ast }(a^{\ast \ast })\right\rangle $ since $i^{\ast \ast }\left[
(B_{0}^{\ast },B_{1}^{\ast })_{\alpha }^{\ast }\cap (B_{0}^{\ast
},B_{1}^{\ast })_{\beta }^{\ast }\right] $ is dense in $B_{\theta ,p}^{\ast
\ast }$, $x_{n}^{\ast }\rightarrow _{n\rightarrow +\infty }x^{\ast }$ in the
weak topology of $B_{\theta ,p}^{\ast }.\blacksquare $

\end{document}